\documentclass[a4paper]{coupled2023}

\usepackage{graphicx}
\usepackage{amsmath}
\usepackage{amsfonts}
\usepackage{amssymb}
\usepackage{amsmath,bm}
\usepackage{xcolor}
\usepackage{url}
\usepackage{cleveref}
\usepackage{todonotes}
\usepackage{pdflscape}
\usepackage[percent]{overpic}    
\usepackage{ifthen}               
\newboolean{boldeqnitalic}
\setboolean{boldeqnitalic}{false}
\usepackage{defdb1}               

\title{A Comparison Of Direct Solvers In FROSch Applied To Chemo-Mechanics}

\author{Alexander Heinlein$^{*}$, Bjoern Kiefer$^{\dag}$, Stefan Prüger$^{\dag}$, Oliver Rheinbach$^{\diamondsuit\vartriangle}$ and Friederike Röver$^{\vartriangle}$}

\heading{A. Heinlein, B. Kiefer, S. Prüger, O. Rheinbach, and F. Röver}

\address{$^{*}$Delft University of Technology\\
	Delft Institute of Applied Mathematics\\
	Mekelweg 4, 2628 CD Delft, Netherlands\\
	e-mail: \url{a.heinlein@tudelft.nl}
	\and
	$^{\dag}$Technische Universität Bergakademie Freiberg\\
	Institute of Mechanics and Fluid Dynamics\\
	Lampadiusstr. 4, 09599 Freiberg, Germany\\
	e-mail: \url{{Bjoern.Kiefer,Stefan.Prueger}@imfd.tu-freiberg.de}
	\and
	$^{\diamondsuit}$Technische Universität Bergakademie Freiberg\\
	Institut für Numerische Mathematik und Optimierung; also
	Center for Efficient High Temperature Processes and Materials Conversion (ZeHS), and University Computing Center (URZ)
	\\
	Akademiestr.6, 09599 Freiberg, Germany\\
	e-mail: \url{Oliver.Rheinbach@math.tu-freiberg.de}
	\and
	$^{\vartriangle}$Technische Universität Bergakademie Freiberg\\
	University Computing Center (URZ)\\
	Bernhard-v.-Cotta Str. 1, 09599 Freiberg, Germany\\
	e-mail: \url{Friederike.Roever@hrz.tu-freiberg.de}
}


\abstract{
	Sparse direct linear solvers are at the computational core of domain decomposition preconditioners and therefore have a strong impact on their performance. In this paper, we consider  the Fast and Robust Overlapping Schwarz (FROSch) solver framework of the Trilinos software library, which contains a parallel implementations of the GDSW domain decomposition preconditioner. We compare three different sparse direct solvers used to solve the subdomain problems in FROSch. The preconditioner is applied to different model problems; linear elasticity and more complex fully-coupled deformation diffusion-boundary value problems from chemo-mechanics. We employ FROSch in fully algebraic mode, and therefore, we do not expect numerical scalability. Strong scalability is studied from 64 to 4\,096 cores, where good scaling results are obtained up to 1\,728 cores. The increasing size of the coarse problem increases the solution time for all sparse direct solvers.
}

\begin{document}
	\thispagestyle{empty}
	
	\section{\texttt{FROSch} preconditioner framework}
	Domain decomposition preconditioners are suitable for parallel computations, since they decompose, based on the computational domain, the problem into smaller subdomain problems, which can be solved in parallel. In this paper, we consider the Fast and Robust Overlapping Schwarz (\texttt{FROSch}) preconditioner framework~\cite{Heinlein:2018:FPI} of the \texttt{Trilinos} software library~\cite{trilinosrepo}. The framework contains a parallel implementation of the Generalized-Dryja-Smith-Widlund (GDSW) preconditioner~\cite{Dohrmann:2008:DDL}. The GDSW preconditioner is a two-level overlapping Schwarz preconditioner~\cite{Toselli:2005:DDM} with an energy-minimizing coarse space, which can be written in the form
	\begin{equation}
		M_{\rm  GDSW}^{-1} ={\Phi K_0^{-1} \Phi^T}+\sum\nolimits_{i = 1}^N R_i^T K_i^{-1} R_i\,.
		\label{eq:gdsw}
	\end{equation}
	Here, $K_i = R_i K R_i^T\,\,, i=1,\ldots, N$, represent the local subdomain problems on the overlapping subdomain $\Omega_i'$, which we solve using a sparse direct solver. Each overlapping subdomain $\Omega_i'$ is obtained by recursively adding layers of elements to the nonoverlapping subdomain $\Omega_i$. The global coarse problem $K_0 = \Phi^T K \Phi$ is solved using a sparse direct solver as well. The matrix $\Phi$ contains the coarse basis functions, as columns, spanning the global coarse space. For the construction of the GDSW coarse space functions, we consider interface functions $\Phi_{\Gamma}$ of the nonoverlapping decomposition. The interface functions are chosen as restrictions of the null space of the global Neumann matrix to the interface components, such as the vertices, edges and in 3D faces.  We obtain the global coarse basis functions $\Phi$ by energy minimizing extensions of $\Phi_{\Gamma}$ into the interior of the nonoverlapping subdomain $\Omega_i$. We  obtain
	\begin{equation} 
		\label{eq:phi}
		\Phi
		=
		\begin{bmatrix}
			-K_{II}^{-1}K_{I\Gamma}\Phi_\Gamma \\
			\Phi_\Gamma
		\end{bmatrix}.
	\end{equation}
	For scalar elliptic problems and regular decompositions the GDSW preconditioner has a known condition number bound
	\begin{equation}
		\kappa (M_{\rm GDSW}^{-1} K) \le C \left(1+ \frac{H}{\delta}\right)\left(1+ \log\left(\frac{H}{h}\right)\right),
		\label{eq:cond}
	\end{equation}
	where $C$ is a constant independent of the problem parameters, $h$ the element, $H$ the subdomain size and $\delta$ the subdomain overlap; cf.~\cite{Dohrmann:2008:FEM,Dohrmann:2008:DDL}. For extended parallel scalability, the \texttt{FROSch} framework includes an implementation of the reduced dimensional GDSW (RGDSW) coarse space~\cite{Dohrmann:2017:DSC} as well as a multi-level extension~\cite{Heinlein:2022:PST}. However, for all results presented here, we only applied the classical GDSW coarse space and two levels. 
	
	\section{Model problems}
	\subsection{Linear elasticity}
	\label{sec:LinearElasticity}
	As a first model problem, we consider the linear elasticity model problem in three dimensions: find ${\bf u} \in H^1(\Omega)^3$
	\begin{align}
		-\div(\boldsymbol{\sigma}) &= f \quad \mbox{in}\,\Omega \\ \nonumber
		\boldsymbol{u} &= 0 \quad \mbox{on} \, \partial\Omega_D
	\end{align}
	We use a generic right-hand-side vector of ones $\begin{pmatrix}
		1 \ldots 1
	\end{pmatrix}^T$ and use the standard implementations from~\cite{PAMM2019} based on~\cite{deal:step-8}.
	
	\subsection{Coupled mechanics-diffusion problems} \label{sec:coupled-mech-diffusion}
	\quad In contrast to the linear boundary value problem introduced in~\cref{sec:LinearElasticity}, the model considered in this section incorporates material and geometrically nonlinear effects in a fully coupled formulation of the mechanical balance of momentum and Fickian diffusion.
	The model employs the rate-type potential from~\cite{BogNatMie:2017:min},
	\begin{equation}
		\label{eq:rate-potential-coupled-mech-diffusion}
		\Pi\left(\dot{\Bvarphi},\dot{v},\bm{J}_{v}\right)=\dd{}{t} \underbrace{\int\nolimits_{\Br}\widehat{\psi}\left(\nabla\Bvarphi,v\right)\mathrm{d}V}_{E\left(\Bvarphi,v\right)} + \underbrace{\int\nolimits_{\Br}\widehat{\phi}\left(\bm{J}_{v};\nabla\Bvarphi,v\right)\mathrm{d}V}_{D\left(\bm{J}_{v}\right)} - P_{\mathrm{ext}}^{\Bvarphi}\left(\dot{\Bvarphi}\right) - P_{\mathrm{ext}}^{\mu}\left(\bm{J}_{v}\right) \ \mbox{,}
	\end{equation}
	in which $E$ denotes the stored energy functional, depending on the deformation $\Bvarphi$ through its material gradient and the swelling volume fraction $v$.
	Furthermore, the dissipation potential functional $D$ associated with the body $\Br$ is a function of the fluid flux $\bm{J}_{v}$ and is additionally parameterized by means of the deformation gradient and the swelling volume fraction.
	The mechanical part of the external load functional, associated volumetric body forces and prescribed tractions, is abbreviated with $P_{\mathrm{ext}}^{\Bvarphi}$, whereas $P_{\mathrm{ext}}^{\mu}$ expresses the corresponding diffusion part, which depends upon the normal component of fluid flux, i.e. $H_{v}=\bm{J}_{v}\cdot\bN$.
	By incorporating the balance of solute volume
	\begin{equation}
		\label{eq:balance-solute-volume}
		\dot{v}=-\Divb{\bm{J}_{v}}  
	\end{equation}
	in~\eqref{eq:rate-potential-coupled-mech-diffusion}, the primary fields can be computed from the two-field minimization principle
	\begin{equation}
		\label{eq:two-field-minimization-mech-diffusion}
		\left\{\dot{\Bvarphi},\bm{J}_{v}\right\} = \mathrm{Arg}\left\{\underset{\dot{\Bvarphi}\in \mathcal{W}_{\dot{\Bvarphi}}}{\mathrm{inf}}\ \underset{\bm{J}_{v}\in \mathcal{W}_{\bm{J}_{v}}}{\mathrm{inf}} \Pi\left(\dot{\Bvarphi},\bm{J}_{v}\right) \right\} \ \mbox{,}
	\end{equation}
	in which the following admissible function spaces are chosen:
	\begin{equation} \label{eq:function-spaces-deformation-flux}
		\begin{aligned}
			\mathcal{W}_{\dot{\Bvarphi}}&=\left\{\dot{\Bvarphi}\in H^{1}\!\left(\Br\right) \vert \ \dot{\Bvarphi}=\dot{\bar{\Bvarphi}} \ \text{on} \ \partial \Br^{\Bvarphi}\right\}, \\ \mathcal{W}_{\bm{J}_{v}}&=\left\{\bm{J}_{v}\in H\left(\Div,\Br\right) \vert \ \bm{J}_{v}\!\cdot\bN=H_{v} \ \text{on} \ \partial \Br^{H_{v}}\right\}.
		\end{aligned}
	\end{equation}
	In the implementation of this model, a free-energy function $\psi$ of Neo-Hookean type in connection with a Flory-Rehner type energy that accounts for the energy due to changes in the swelling volume fraction, and a quadratic dissipation potential $\phi$ are chosen.
	Upon the application of the Euler backward time integration to~\eqref{eq:balance-solute-volume} and \eqref{eq:rate-potential-coupled-mech-diffusion}, the time discrete counterpart of~\eqref{eq:two-field-minimization-mech-diffusion} is employed to compute the primary fields $\Bvarphi\vert_{n+1}$ and $\bm{J}_{v}\vert_{n+1}$ at time $t_{n+1}$.
	Note that the employed variational principle ensures that the linearization of the necessary optimality condition yields a symmetric system of equations.
	For more details, the interested reader is referred to~\cite{Kiefer:2023:MPO,BogNatMie:2017:min}. The following specific model problems were also used in~\cite{Kiefer:2023:MPO}.
	
	\begin{figure}[t]
		\centering
		\begin{overpic}[scale=1.1]{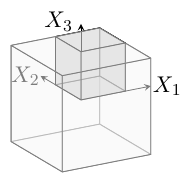}
			\put(0,0){(a)}
		\end{overpic}
		\hspace{20ex}
		\begin{overpic}[scale=1.1]{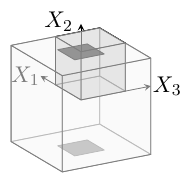}
			\put(0,0){(b)}
		\end{overpic}
		\caption{Cuboidal domain considered in the free-swelling boundary problem (a), and mechanically induced diffusion boundary value problem. Figures taken from~\cite[Figure 4 and 6]{Kiefer:2023:MPO}.}
		\label{fig:boundary-value-problems}
	\end{figure}
	
	\subsubsection{Free-swelling boundary value problem} \label{sec:free-swelling-bvp}
	This problem was also investigated in~\cite{Kiefer:2023:MPO}. In the free-swelling boundary value problem, a cube of edge length $2L$, as shown in~\Cref{fig:boundary-value-problems}(a), is considered. It is loaded in terms of a temporarily varying fluid flux at the outer boundary, while the outer surface remains traction free.
	Due to the intrinsic symmetry of the problem, only one eighth of the cube is taken into account in the scalability studies.
	Therefore, appropriate symmetry boundary conditions are applied along the symmetry planes ($X_{1}=0, X_{2}=0$, and $X_{3}=0$), i.e., the displacement component and the fluid flux in the direction normal to these planes are set to zero.
	
	\subsubsection{Mechanically induced diffusion boundary value problem} \label{sec:mech-induced-bvp}
	This problem was also investigated in~\cite{Kiefer:2023:MPO}. Similar to the free-swelling boundary value problem, a cuboidal domain is also considered for the mechanically induced diffusion problem. However, here, zero Neumann boundary conditions for the normal component of the fluid flux are prescribed, while along the subset $\left(X_{1},X_{3}\right)\in\left[-\frac{\mathrm{L}}{3},\frac{\mathrm{L}}{3}\right]\times\left[-\frac{\mathrm{L}}{3},\frac{\mathrm{L}}{3}\right]$ at the plane $X_{2}=\mathrm{L}$, highlighted in~\Cref{fig:boundary-value-problems}(b), the coefficients of the displacement vector are prescribed as $u_{i}=[0,-\hat{u},0]$. Once again the intrinsic symmetry of the problem is exploited by specifying symmetry boundary conditions as described in~\cref{sec:free-swelling-bvp}. The material parameters employed in the free energy function $\psi$ and the dissipation potential $\phi$ are adopted from~\cite{Kiefer:2023:MPO}.
	
	\section{Implementation}
	\quad In this paper, we use the implementation and setup from~\cite{Kiefer:2023:MPO}.
	For the the finite element implementation of our model problem, we employed the \texttt{deal.II} software library~\cite{dealiiParallel} version 9.2.0. The MPI-parallel data distribution is handled using the \texttt{parallel::distributed::Triangulation}, which links to the external \texttt{p4est} library~\cite{p4est}. To solve the nonlinear system, we apply the Newton-Raphson scheme with the relative and absolute tolerances as provided in~\Cref{tab:tol-newton}.
	\begin{table}[t]
		\centering
		\begin{tabular}{c|c|c}
			& $\Vert r_k \Vert$ & $\Vert r_k \Vert/\Vert r_0 \Vert$ \\ \hline
			$\boldsymbol{\varphi}$ & $10^{-9}$ & $10^{-6}$ \\
			$\bm{J}_{v}$ & $5 *10 ^{-12}$ & $10^{-9}$ \\ 
		\end{tabular}
		\caption{Tolerances for the Newton-Raphson scheme. Here, $r_k$ is the $k$-th residual. Table taken from~\cite[Table 3]{Kiefer:2023:MPO}}
		\label{tab:tol-newton}
	\end{table}
	For the parallel linear algebra, we use the \texttt{deal.II} \texttt{TrilinosWrappers} such that the \texttt{Trilinos} package \texttt{Epetra} is applied. Further, we use some functions implemented for standard tensor operations from~\cite{McBJavSteRed:2015:fin}. The linearized system is solved using the Krylov iteration method GMRES. For GMRES, we use the parallel implementation provided by the \texttt{Trilinos} package \texttt{Belos} with a relative stopping criterion of $\Vert r_k \Vert/\Vert r_0 \Vert \leq 10^{-8}$, where $r_k$ is the $k$-th residual and $r_0$ the initial residual. \texttt{FROSch} is applied as a preconditioner. In all computations, we apply an algebraically computed overlap of two nodes ($\delta \approx 2h$) and employ the provided algebraic computation of the interface components. Applying \texttt{FROSch} in a completely algebraic sense implies using a one-dimensional null space. It has been shown that \texttt{FROSch} may be able to scale even if certain dimensions of the null space are neglected~\cite{Heinlein:2016:PIT,heinlein:2021:FAT}. However, this is not covered by the theory. A one-to-one correspondence between cores and subdomains is applied, and the global coarse problem is solved on a single core. We use \texttt{Trilinos} version 13.0.1 with small modifications. We compare the performance of different sparse direct linear solvers, applying the build-in \texttt{KLU} solver from \texttt{Amesos2} as well as \texttt{Umfpack}~\cite{umfpack} and \texttt{MUMPS}~\cite{mumps} both interfaced through \texttt{Amesos}. We always use the same sparse solver for the local problems and the coarse problem. Using \texttt{Trilinos} version 13.0.1, we faced issues with the \texttt{Amesos2} interfaces for \texttt{Epetra} matrices, such that we used the older \texttt{Amesos} interfaces in our tests; in more recent \texttt{Trilinos} versions, these issues may have been fixed. We use \texttt{MUMPS} in Version 5.6.0 without \texttt{METIS}, and \texttt{Umfpack} included in the \texttt{Suite Sparse} library Version 5.1.0, which uses \texttt{METIS} 5.1.0-IDX64. We consider the \textit{solver time}, which denotes the time to build the preconditioner and to perform the Krylov iterations. The time to solve the subdomain problems as well as the time to solve the coarse problem includes the numerical factorization and the forward/backward substitution, denoted as \textit{subd.~problems solve time} and \textit{coarse problem solve time}, respectively.  For the \textit{coarse problem solve time} the time is determined by lower level timers, such that it may deviate from the pure solution time by the sparse direct linear solver. All test were performed on the JSC supercomputer JUWELS~\cite{JUWELS} at the FZ J\"ulich using the Intel 2021.4.0 compiler and IntelMPI.
	
	\section{Numerical Results}
	\label{sec:num-results}
	
	\begin{figure}[t]
		\includegraphics[width=0.48\textwidth]{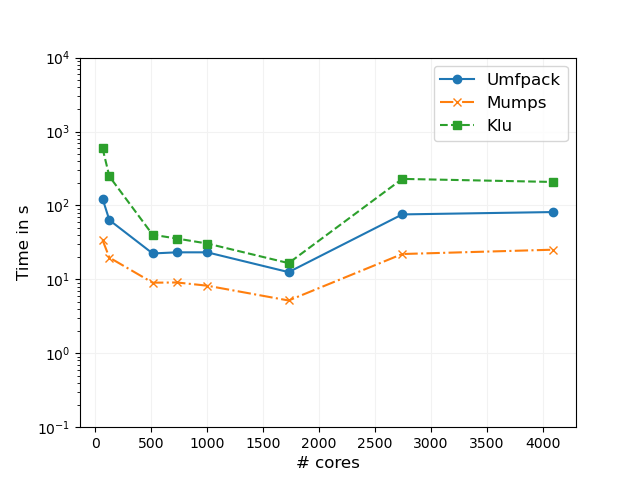}
		\includegraphics[width=0.48\textwidth]{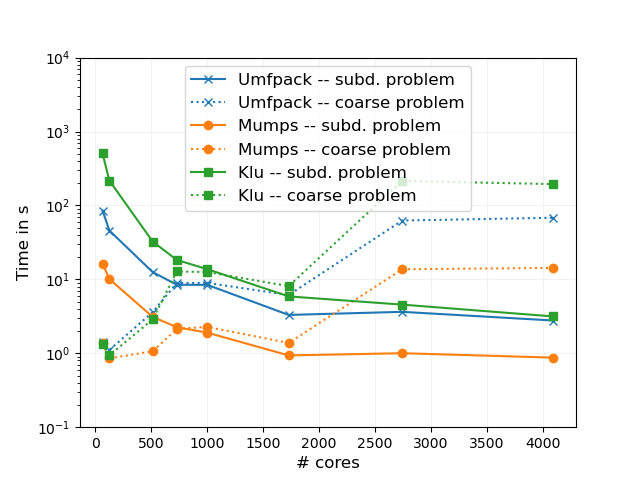}
		\caption{Strong scalability of the \textit{solver time}, scaling from 64 to 4\,096 cores, for the linear elasticity model problem applying the \texttt{FROSch} preconditioner with  th GDSW coarse space (left) and the detailed timers \textit{subd. problem solve time} and the \textit{coarse problem solve time} (right), including the time for the factorization and the forward./backward substitution. See~\Cref{tab:linElas} for the data.
		\label{fig:linEals}}
	\end{figure}
	
	\begin{table}[t]
		\centering
		\renewcommand{\arraystretch}{1.5}
		\fontsize{6pt}{6pt}\selectfont \addtolength{\tabcolsep}{-0.5mm}
		\begin{tabular}{|r|r|r|r||rrr|rrr|rrr|}
			\hline
			\multicolumn{13}{|c|}{\bf linear elasticity -- GDSW}                                                                                                                                                                \\ \hline
			\multicolumn{4}{|c||}{}                  & \multicolumn{3}{|c|}{\bf solver time}                  & \multicolumn{3}{|c|}{\bf subd.~problems}               & \multicolumn{3}{|c|}{\bf coarse problem}               \\
			\multicolumn{4}{|c||}{}                  & \multicolumn{3}{|c|}{}                                 & \multicolumn{3}{|c|}{\bf solve time}                   & \multicolumn{3}{|c|}{\bf solve time}                   \\ \hline
			\# cores & Krylov &       max. &    size & \multicolumn{2}{c}{\texttt{Amesos}} & \texttt{Amesos2} & \multicolumn{2}{c}{\texttt{Amesos}} & \texttt{Amesos2} & \multicolumn{2}{c}{\texttt{Amesos}} & \texttt{Amesos2} \\
			&    it. & size $K_i$ &   $K_0$ & \texttt{MUMPS} &   \texttt{Umfpack} &     \texttt{KLU} & \texttt{MUMPS} &   \texttt{Umfpack} &     \texttt{KLU} & \texttt{MUMPS} &   \texttt{Umfpack} &     \texttt{KLU} \\ \hline\hline
			64 &     56 &    86\,577 &     932 &         34.04s &          123.38\,s &        592.77\,s &       16.36\,s &           85.23\,s &        505.00\,s &        1.36\,s &            1.44\,s &          1.34\,s \\
			125 &     62 &    54\,396 &  2\,108 &         19.72s &           63.39\,s &        246.96\,s &       10.11\,s &           45.80\,s &        217.09\,s &        0.85\,s &            1.11\,s &          0.94\,s \\
			512 &     72 &    21\,504 & 10\,412 &          9.02s &           22.41\,s &         40.16\,s &        3.11\,s &           12.54\,s &         32.50\,s &        1.07\,s &            3.59\,s &          2.87\,s \\
			729 &     73 &    17\,868 & 16\,412 &          9.07s &           23.21\,s &         35.64\,s &        2.25\,s &            8.44\,s &         18.27\,s &        2.14\,s &            8.95\,s &         12.79\,s \\
			1\,000 &     78 &    14\,724 & 20\,037 &          8.21s &           23.21\,s &         30.54\,s &        1.90\,s &            8.44\,s &         13.74\,s &        2.28\,s &            8.95\,s &         12.56\,s \\
			1\,728 &     68 &     7\,581 & 18\,788 &          5.18s &           12.53\,s &         16.58\,s &        0.94\,s &            3.30\,s &          5.88\,s &        1.38\,s &            6.14\,s &          8.10\,s \\
			2\,744 &     86 &     8\,700 & 60\,090 &         22.02s &           75.62\,s &        228.62\,s &        1.00\,s &            3.64\,s &          4.55\,s &       13.69\,s &           62.47\,s &        214.63\,s \\
			4\,096 &     90 &     7\,038 & 78\,653 &         25.20s &           81.41\,s &        207.76\,s &        0.87\,s &            2.77\,s &          3.15\,s &       14.29\,s &           68.19\,s &        194.30\,s \\ \hline
		\end{tabular}
		\caption{Strong scalability results for the linear elasticity model problem using $Q_1$ elements. We apply the \texttt{FROSch} preconditioner with the GDSW coarse space and an algebraic overlap of two elements and compare different sparse direct linear solver. The \textit{solver time} is the time to build the preconditioner and to perform the Krylov iterations. The problem size is 2\,738\,019.
			\label{tab:linElas}}
	\end{table}
	
	\subsection{Linear Elasticity}
	As a first example problem, we choose the linear elasticity model problem, described in~\cref{sec:LinearElasticity}, using $Q_1$ finite elements on a structured mesh with 884\,736 cells such that we have 2\,738\,019 degrees of freedom (DOF). Since we neglect the rotations from the null space, we may expect the number of iterations to increase with an increasing number of cores. For our tests, the number of iterations increases from 56 to 90 scaling from 64 to 4\,096 cores; see~\Cref{tab:linElas}. Note that it has shown that the majority of the \textit{solver time} is taken by the construction of the preconditioner rather than by the Krylov iterations for a smaller number  of cores~\cite{Kiefer:2023:MPO}. For this reason, the number of iterations should not influence the \textit{solver time} significantly compared to the sizes of the subdomain problems within the considered range of cores.
	
	We obtain good strong scalability results scaling from 64 to 1\,728 cores for all sparse direct linear solvers with \texttt{MUMPS} being the fastest. \texttt{Umfpack} is faster than \texttt{KLU}, however for smaller subdomain problem sizes (obtained from 729 to 1\,728 cores) the results are comparable; see~\Cref{tab:linElas} and~\Cref{fig:linEals}. Reaching 2\,744 cores the \textit{solver time} starts to increase instead of decreasing. The reason for that is the 
	significant increase for the size of the coarse problem, e.g., for 1\,728 dim($K_0$) =18\,788 which compares to dim($K_0$) = 60\,090 for 2\,744 cores; see~\Cref{tab:linElas}. From 2\,744 cores on, the problem size relation of $K_0$ and max $K_i$ is similar to 125 cores, where the \textit{subd. problem solve time} dominated the \textit{solver time}. The coarse problem reaches a critical size beyond 1\,728 cores, since after this point the solution of the coarse problem begins to dominate the \textit{solver time}; see~\Cref{fig:linEals}. The solution on a single core is not sufficient anymore. As a next step, we could apply three-levels or/and RGDSW. We want to remark that for certain numbers of cores and numbers of elements the decomposition obtained from \textit{p4est} results in structured decompositions decreasing the size of the coarse problem. The proportion of the number of elements and the number of subdomains is decisive for this phenomenon.
	
	\begin{table}[t]
		\centering
		\renewcommand{\arraystretch}{1.5} 
		\fontsize{6pt}{6pt}\selectfont \addtolength{\tabcolsep}{-0.5mm}
		\begin{tabular}{|r||r|rrr||r|rrr|}
			\hline
			\multicolumn{9}{|c|}{\bf free-swelling problem -- GDSW}                                                                                            \\
			\multicolumn{9}{|c|}{\bf solver time}                                                                                                          \\ \hline
			&                  \multicolumn{4}{c||}{$Q_1Q_1$}                  & \multicolumn{4}{|c|}{$Q_1RT_0$}                                  \\ \hline
			&   avg. & \multicolumn{2}{|c}{\texttt{Amesos}} & \texttt{Amesos2} &   avg. & \multicolumn{2}{|c}{\texttt{Amesos}} & \texttt{Amesos2} \\
			\# cores & Krylov & \texttt{MUMPS} &    \texttt{Umfpack} &     \texttt{KLU} & Krylov & \texttt{MUMPS} &    \texttt{Umfpack} &     \texttt{KLU} \\ \hline
			64 &   70.9 &      173.24\,s &           536.41\,s &     1\,587.04\,s &   49.7 &       75.35\,s &           238.57\,s &        686.13\,s \\
			125 &   79.8 &      121.68\,s &           362.16\,s &        833.82\,s &   56.4 &       53.18\,s &           155.20\,s &        352.30\,s \\
			512 &  103.4 &       71.88\,s &           213.51\,s &        250.88\,s &   73.7 &       41.18\,s &           122.68\,s &        134.44\,s \\
			729 &  111.0 &       78.86\,s &           258.50\,s &        240.84\,s &   78.3 &       54.20\,s &           178.32\,s &        200.62\,s \\
			1\,000 &  116.3 &       73.49\,s &           241.98\,s &        216.39\,s &   81.3 &       53.48\,s &           167.36\,s &        154.35\,s \\
			1\,728 &  105.1 &       57.20\,s &           190.70\,s &        171.00\,s &   76.5 &       38.76\,s &           121.18\,s &        128.17\,s \\
			2\,744 &  125.1 &      147.99\,s &           599.51\,s &     1\,068.12\,s &   86.8 &      172.00\,s &           646.72\,s &     1\,547.63\,s \\
			4\,096 &  131.9 &      178.61\,s &           755.28\,s &     1\,428.17\,s &   92.3 &      197.48\,s &           836.19\,s &     1\,428.17\,s \\ \hline
		\end{tabular}
		\caption{Strong scalability results for the free-swelling model problem. We apply the \texttt{FROSch} preconditioner with the GDSW coarse space and an algebraic overlap of two elements and compare different sparse direct linear solver. The \textit{solver time} is the time to build the preconditioner and to perform the Krylov iterations. The problem size is 705\,894 for $Q_1Q_1$ elements and 691\,635 for $Q_1RT_0$ elements.}
		\label{tab:free-solver}
	\end{table}

	\begin{figure}[t]
		\includegraphics[width=0.48\textwidth]{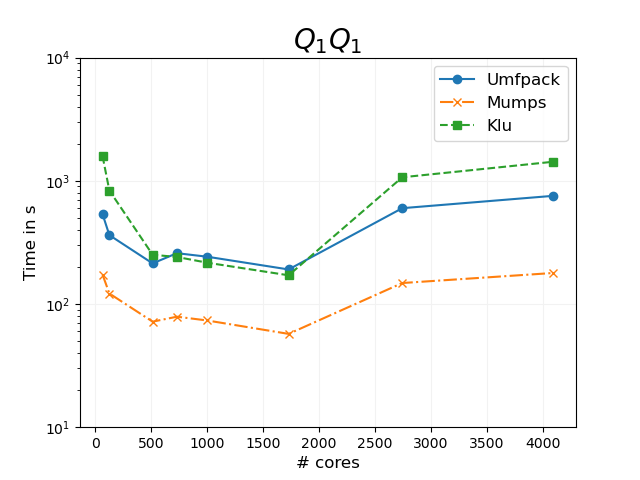}
		\includegraphics[width=0.48\textwidth]{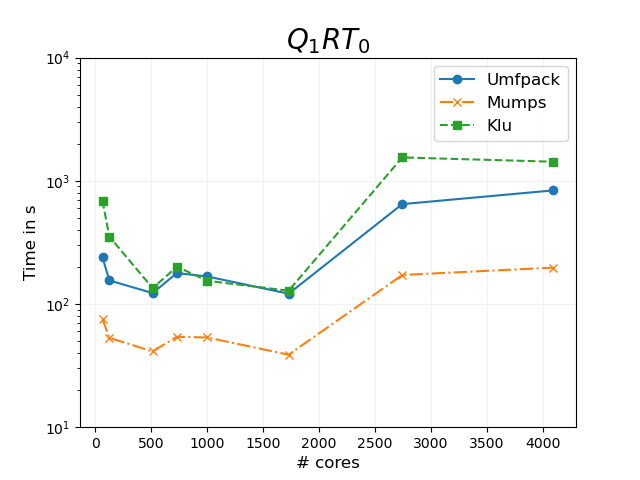}
		\caption{\textit{solver time} scaling from 64 to 4\,096 cores using different direct linear solvers for the  subd.~problems of the \texttt{FROSch} preconditioner with a GDSW coarse space. The preconditioner is applied the \textit{free-swelling} model problem. See~\Cref{tab:free-solver} for the data.
			\label{fig:free-compare-solver}}
	\end{figure}
	
	\subsection{Free-swelling boundary value problem}
	For the model problem described in~\cref{sec:free-swelling-bvp}, we consider a mesh with 110\,592 cells. We compare two types of ansatz functions for the flux field $Q_1$ and $RT_0$, and  for the deformation field, we always use $Q_1$ elements. For $Q_1Q_1$, we obtain 705\,894 DOF and 691 635 DOF for  $Q_1RT_0$. We restrict the computation to two time steps. Each requires five Newton iterations. We take the \textit{solver time} over the two time steps, such that the preconditioner is applied ten times. We consider an average number of Krylov iterations (\textit{avg.~Krylov}) over the ten Newton steps. Since we apply \texttt{FROSch} fully algebraically, i.e., using the the null space of the Laplace operator for the construction of the coarse space, we cannot expect numerical scalability. Consequently, \textit{avg.~Krylov} increases with the number of cores; see~\Cref{tab:free-solver}.
	
	Generally, the number of Krylov iterations is lower for the $Q_1RT_0$ elements, yet the increase from 64 to 4\,096 cores stronger than for the $Q_1Q_1$ elements. This is in agreement with the results obtained in~\cite{Kiefer:2023:MPO}. Further, the size of the coarse problem is much larger with increasing number of cores. For 4\,096 cores, we obtain a coarse space dimension of 48\,232 for $Q_1Q_1$ elements respectively 77\,222 for $Q_1RT_0$ elements; see also~\Cref{tab:free-detail}. This indicates that the algebraic decomposition for the $Q_1Q_1$ is favorable since less interface components are obtained.
	
\begin{figure}[t]
		\includegraphics[width=0.48\textwidth]{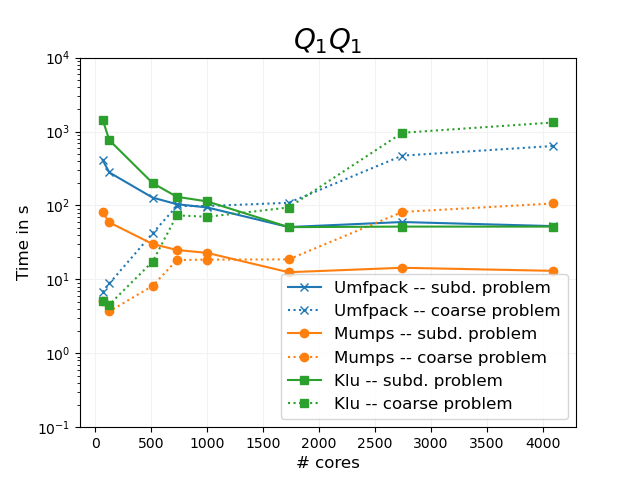}
		\includegraphics[width=0.48\textwidth]{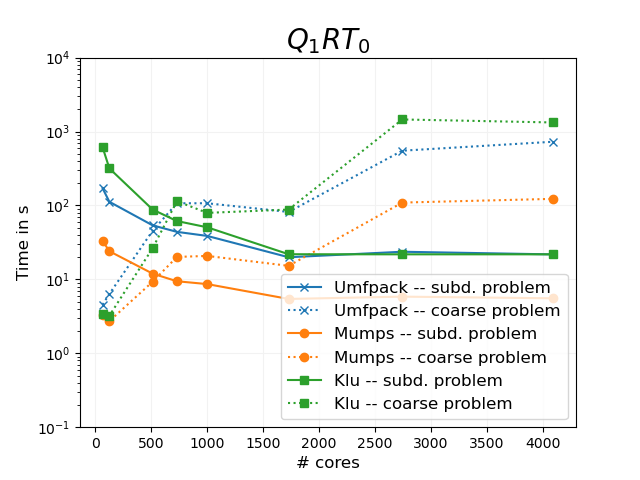}
		\caption{Time to solve the local subd.~problems scaling from 64 to 4\,096 cores using different direct linear. See \Cref{tab:free-detail} for the data.}
		\label{fig:free-compare-osolver}
	\end{figure}

	Regarding the strong scalability, we obtain similar results as for the linear elasticity model problem; compare~\Cref{fig:linEals,fig:free-compare-solver}. For both ansatz functions, the \textit{solver time} increases reaching 2\,744 cores. The $Q_1RT_0$ elements are faster up to 1\,728 cores. For larger number of cores the \textit{solve coarse problem time} dominates the \textit{solver time}. As previously discussed, the size of the coarse problem is smaller for $Q_1Q_1$ elements such that the \textit{solver time} is faster for these elements employing larger numbers of cores. As for linear elasticity, the best performance of the solver framework is obtained using \texttt{MUMPS}. For 64 cores, \texttt{MUMPS} is more than 15 times faster than \texttt{KLU} and more than five times faster than \texttt{Umfpack}; see~\Cref{tab:free-solver}. The advantage of using \texttt{MUMPS} is most apparent if the size of the directly solved problem is large. From 1\,728 to 4\,096 the \textit{subd.~problem solve time} \texttt{Umfpack} and \texttt{KLU} perform similarly, e.g., considering 1\,728 cores, the time to solve the $K_i$s is 51.02s for \texttt{Umfpack} and 50.89s for \texttt{KLU} using $Q_1Q_1$ elements; see also~\Cref{tab:free-detail} and~\Cref{fig:free-compare-osolver}. Generally for smaller subdomain problem sizes \texttt{Umfpack} and \texttt{KLU} are comparable.
	
	\begin{table}[t]
		\centering
		\renewcommand{\arraystretch}{1.5}
		\fontsize{6pt}{6pt}\selectfont \addtolength{\tabcolsep}{-0.5mm}
		{
			\begin{tabular}{|r||r|rrr||r|rrr|}
				\hline
				\multicolumn{9}{|c|}{\bf free-swelling problem -- GDSW}                                                                                            \\ \hline
				&                   \multicolumn{4}{c||}{$Q_1Q_1$}                    & \multicolumn{4}{|c|}{$Q_1RT_0$}                                     \\ \hline
				\multicolumn{9}{|c|}{\bf subd.~problem solve time}                                                                                                 \\ \hline
				& max.~size & \multicolumn{2}{|c}{\texttt{Amesos}} & \texttt{Amesos2} & max.~size & \multicolumn{2}{|c}{\texttt{Amesos}} & \texttt{Amesos2} \\
				&     $K_i$ & \texttt{MUMPS} &    \texttt{Umfpack} &     \texttt{KLU} &     $K_i$ & \texttt{MUMPS} &    \texttt{Umfpack} &     \texttt{KLU} \\ \hline
				64 &   38\,334 &       81.29\,s &           412.17\,s &     1\,453.55\,s &   35\,634 &       33.33\,s &           171.26\,s &        615.22\,s \\
				125 &   27\,648 &       58.90\,s &           279.19\,s &        763.12\,s &   25\,598 &       24.02\,s &           112.08\,s &        317.11\,s \\
				512 &   13\,896 &       30.06\,s &           127.59\,s &        200.50\,s &   12\,571 &       11.95\,s &            53.77\,s &         87.94\,s \\
				729 &   11\,472 &       24.97\,s &           103.84\,s &        130.83\,s &   10\,368 &        9.43\,s &            43.92\,s &         61.50\,s \\
				1000 &   10\,452 &       22.79\,s &            93.92\,s &        113.40\,s &    9\,369 &        8.63\,s &            38.68\,s &         50.62\,s \\
				1728 &    5\,634 &       12.49\,s &            51.02\,s &         50.89\,s &    5\,064 &        5.43\,s &            19.91\,s &         21.87\,s \\
				2744 &    7\,020 &       14.35\,s &            59.66\,s &         51.74\,s &    6\,209 &        5.84\,s &            23.64\,s &         21.81\,s \\
				4096 &    6\,312 &       13.03\,s &            52.51\,s &         51.74\,s &    5\,563 &        5.54\,s &            21.75\,s &         21.81\,s \\ \hline
				\multicolumn{9}{|c|}{\bf coarse problem solve time}                                                                                                \\ \hline
				&      size & \multicolumn{2}{|c}{\texttt{Amesos}} & \texttt{Amesos2} &      Size & \multicolumn{2}{|c}{\texttt{Amesos}} & \texttt{Amesos2} \\
				&     $K_0$ & \texttt{MUMPS} &    \texttt{Umfpack} &     \texttt{KLU} &     $K_0$ & \texttt{MUMPS} &    \texttt{Umfpack} &     \texttt{KLU} \\ \hline
				64 &       981 &        5.07\,s &             6.74\,s &          5.16\,s &       981 &        3.35\,s &             4.53\,s &          3.39\,s \\
				125 &    2\,003 &        3.69\,s &             8.80\,s &          4.56\,s &    2\,125 &        2.72\,s &             6.35\,s &          3.22\,s \\
				512 &    8\,614 &        8.16\,s &            42.14\,s &         17.41\,s &   10\,748 &        9.30\,s &            45.55\,s &         26.41\,s \\
				729 &   13\,552 &       18.24\,s &            98.23\,s &         73.48\,s &   16\,712 &       20.11\,s &           107.02\,s &        114.12\,s \\
				1\,000 &   15\,904 &       18.51\,s &            98.23\,s &         70.23\,s &   20\,578 &       20.76\,s &           107.02\,s &         79.36\,s \\
				1\,728 &   18\,788 &       18.62\,s &           108.68\,s &         93.77\,s &   18\,788 &       15.24\,s &            81.61\,s &         87.93\,s \\
				2\,744 &   37\,229 &       81.63\,s &           471.65\,s &        961.46\,s &   57\,308 &      109.34\,s &           550.39\,s &     1\,456.04\,s \\
				4\,096 &   48\,232 &      106.39\,s &           637.60\,s &     1\,325.53\,s &   77\,222 &      122.81\,s &           727.94\,s &     1\,325.53\,s \\ \hline
			\end{tabular}
		}
		\caption{Strong scalability for the \textit{subd. problem solve time} and the \textit{coarse problem solve time} using different sparse direct linear solvers. The time includes the factorization of problem and the forward./backw. substitution during the Krylov iterations. For the whole \textit{solver time} see~\Cref{tab:free-solver}}
		\label{tab:free-detail}
	\end{table}
	
	Although \texttt{MUMPS} is much faster than \texttt{Umfpack} and \texttt{KLU} the scalability is not extended solely by this choice of the direct linear solver.

	\subsection{Mechanically induced diffusion problem} \label{sec:mechnially-induced}
\begin{figure}[t]
		\includegraphics[width=0.48\textwidth]{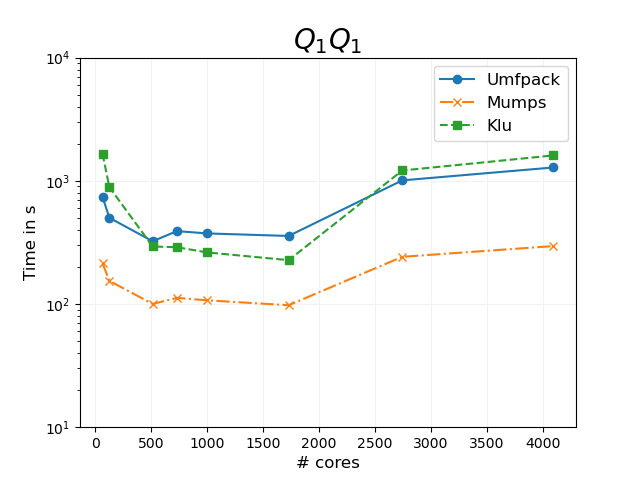}
		\includegraphics[width=0.48\textwidth]{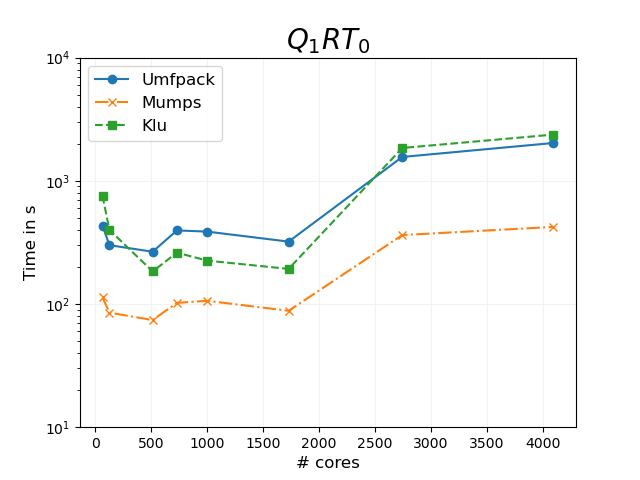}
		\caption{\textit{Solver time} scaling from 64 to 4\,096 cores using different direct linear solvers for the  subd.~problems of the \texttt{FROSch} preconditioner with a GDSW coarse space. The preconditioner is applied the \textit{mechanically induced diffusion} model problem. See~\Cref{tab:mech-solver} for the data.}
		\label{fig:mech-compare-solver}
	\end{figure}

\begin{figure}[t]
		\includegraphics[width=0.48\textwidth]{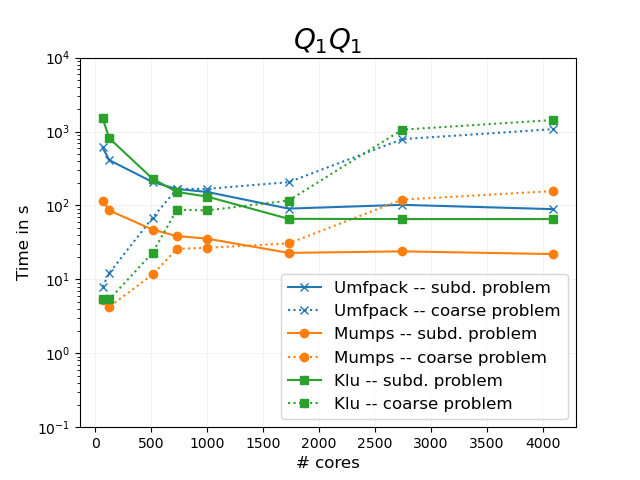}
		\includegraphics[width=0.48\textwidth]{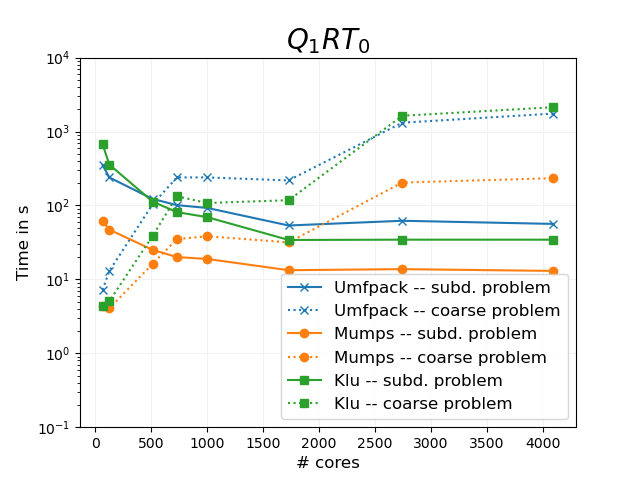}
		\caption{Time to solve the local subd.~problems scaling from 64 to 4\,096 cores using different direct linear. See~\Cref{tab:mech-detail} for the data.}
		\label{fig:mech-compare-osolver}
	\end{figure}

	The results for the mechanically induced diffusion problem, introduced in~\cref{sec:mech-induced-bvp}, confirm the results obtained for the free-swelling boundary value problem. Here, we also restrict the computation to two timesteps each solved with five Newton iterations.
	The more complex boundary conditions result in higher numbers of Krylov iterations. As for the other model problems, the number os iterations increases with an increasing number of cores.
	\begin{table}[t]
		\centering
		\renewcommand{\arraystretch}{1.5}
		\fontsize{6pt}{6pt}\selectfont \addtolength{\tabcolsep}{-0.5mm}
		\begin{tabular}{|r||r|rrr||r|rrr|}
			\hline
			\multicolumn{9}{|c|}{\bf mechanically-induced diffusion problem -- GDSW}                                                  \\ 
			\multicolumn{9}{|c|}{\bf solver time}                                                                               \\ \hline
			&           \multicolumn{4}{c||}{$Q_1Q_1$}            & \multicolumn{4}{|c|}{$Q_1RT_0$}                     \\ \hline
			&   avg. & \multicolumn{2}{|c}{\texttt{Amesos}} &      \texttt{Amesos2} &   Avg. & \multicolumn{2}{|c}{\texttt{Amesos}} &      \texttt{Amesos2} \\
			\# cores & Krylov &     \texttt{MUMPS} &         \texttt{Umfpack} &          \texttt{KLU} & Krylov &     \texttt{MUMPS} &         \texttt{Umfpack} &          \texttt{KLU} \\ \hline
			64 &  115.8 & 214.14\,s &       739.30\,s & 1\,654.99\,s &  119.9 & 113.30\,s &       432.78\,s &    753.36\,s \\
			125 &  130.0 & 154.23\,s &       500.88\,s &    887.11\,s &  137.2 &  85.14\,s &       300.36\,s &    399.83\,s \\
			512 &  177.2 & 100.06\,s &       323.04\,s &    292.94\,s &  180.3 &  74.06\,s &       265.61\,s &    183.75\,s \\
			729 &  186.4 & 112.14\,s &       390.50\,s &    288.87\,s &  192.3 & 102.14\,s &       395.56\,s &    260.95\,s \\
			1\,000 &  196.4 & 107.12\,s &       374.02\,s &    262.23\,s &  205.4 & 106.15\,s &       386.94\,s &    224.96\,s \\
			1\,728 &  210.6 &  97.71\,s &       356.62\,s &    226.93\,s &  219.5 &  88.08\,s &       320.87\,s &    192.73\,s \\
			2\,744 &  224.4 & 241.59\,s &    1\,008.08\,s & 1\,212.64\,s &  235.6 & 362.65\,s &    1\,561.48\,s & 1\,846.48\,s \\
			4\,096 &  239.1 & 295.22\,s &    1\,283.16\,s & 1\,606.95\,s &  249.0 & 422.36\,s &    2\,027.79\,s & 2\,372.35\,s \\ \hline
		\end{tabular}
		\caption{Strong scalability results for the \textit{mechanically-induced diffusion} model problem. We apply the \texttt{FROSch} preconditioner with the GDSW coarse space and an algebraic overlap of two elements and compare different sparse direct linear solver. The \textit{solver time} is the time to build the preconditioner and to perform the Krylov iterations. The problem size is 705\,894 for $Q_1Q_1$ elements and 691\,635 for $Q_1RT_0$ elements.}
		\label{tab:mech-solver}
	\end{table} 
	Yet, the change of boundary conditions does not affect the domain decomposition, such that the sizes for the subdomain problems $K_i$ and the coarse problem $K_0$ are equal. Therefore, the strong scaling behavior is similar to the free-swelling problem; see~\Cref{tab:mech-solver}. Consequently, for these tests, the \textit{solver time} obtained using \texttt{MUMPS} is the fastest as well. The results comparing \texttt{Umfpack} and \texttt{KLU} for this model problem are remarkable. From 512 to 1\,728 cores the \textit{solver time} using \texttt{KLU} is slightly faster than \texttt{Umfpack},  although it has been slower for the \textit{free-swelling model problem}. The time to solve the subdomain problems again dominates \textit{solver time}  scaling from 64 to 1\,000 cores; see~\Cref{tab:mech-solver,tab:mech-detail} and~\Cref{fig:mech-compare-solver,fig:mech-compare-osolver}. 
	\begin{table}[t]
		\centering
		\renewcommand{\arraystretch}{1.5}
		\fontsize{6pt}{6pt}\selectfont \addtolength{\tabcolsep}{-0.5mm}
		{
			\begin{tabular}{|r||r|rrr||r|rrr|}
				\hline
				\multicolumn{9}{|c|}{\bf mechanically induced diffusion problem -- GDSW}                                                             \\ \hline
				&             \multicolumn{4}{c||}{$Q_1Q_1$}             & \multicolumn{4}{|c|}{$Q_1RT_0$}                        \\ \hline
				\multicolumn{9}{|c|}{\bf subd.~problem solve time}                                                                          \\ \hline
				& max.~size & \multicolumn{2}{|c}{\texttt{Amesos}} &      \texttt{Amesos2} & max.~size & \multicolumn{2}{|c}{\texttt{Amesos}} &      \texttt{Amesos2} \\
				&     $K_i$ &     \texttt{MUMPS} &         \texttt{Umfpack} &          \texttt{KLU} &     $K_i$ &     \texttt{MUMPS} &         \texttt{Umfpack} &          \texttt{KLU} \\ \hline
				64 &   38\,334 & 115.89\,s &       609.33\,s & 1\,515.37\,s &   35\,634 &  61.56\,s &       351.75\,s &    673.00\,s \\
				125 &   27\,648 &  85.54\,s &       409.80\,s &    809.69\,s &   25\,598 &  47.02\,s &       239.01\,s &    356.35\,s \\
				512 &   13\,896 &  46.94\,s &       207.11\,s &    229.25\,s &   12\,571 &  25.06\,s &       122.50\,s &    112.68\,s \\
				729 &   11\,472 &  38.56\,s &       166.92\,s &    152.12\,s &   10\,368 &  20.03\,s &       100.85\,s &     81.18\,s \\
				1\,000 &   10\,452 &  35.56\,s &       152.18\,s &    131.74\,s &    9\,369 &  18.88\,s &        92.55\,s &     69.41\,s \\
				1\,728 &    5\,634 &  22.80\,s &        90.63\,s &     65.87\,s &    5\,064 &  13.30\,s &        53.54\,s &     34.12\,s \\
				2\,744 &    7\,020 &  23.92\,s &       101.89\,s &     65.43\,s &    6\,209 &  13.75\,s &        62.04\,s &     34.43\,s \\
				4\,096 &    6\,312 &  21.98\,s &        89.07\,s &     65.43\,s &    5\,563 &  13.01\,s &        56.25\,s &     34.43\,s \\ \hline
				\multicolumn{9}{|c|}{\bf coarse problem solve time}                                                                      \\ \hline
				& max.~size & \multicolumn{2}{|c}{\texttt{Amesos}} &      \texttt{Amesos2} & max.~size & \multicolumn{2}{|c}{\texttt{Amesos}} &      \texttt{Amesos2} \\
				&     $K_i$ &     \texttt{MUMPS} &         \texttt{Umfpack} &          \texttt{KLU} &     $K_i$ &     \texttt{MUMPS} &         \texttt{Umfpack} &          \texttt{KLU} \\ \hline
				64 &       981 &   5.32\,s &         7.92\,s &      5.42\,s &       981 &   4.36\,s &         7.18\,s &      4.43\,s \\
				125 &    2\,003 &   4.25\,s &        12.21\,s &      5.34\,s &    2\,125 &   4.06\,s &        13.04\,s &      5.08\,s \\
				512 &    8\,614 &  11.64\,s &        68.46\,s &     22.99\,s &   10\,748 &  16.23\,s &       103.86\,s &     38.81\,s \\
				729 &   13\,552 &  25.86\,s &       167.40\,s &     87.50\,s &   16\,712 &  35.08\,s &       239.38\,s &    132.98\,s \\
				1\,000 &   15\,904 &  26.73\,s &       167.40\,s &     85.64\,s &   20\,578 &  38.30\,s &       239.38\,s &    107.88\,s \\
				1\,728 &   18\,788 &  30.67\,s &       206.60\,s &    116.33\,s &   18\,788 &  31.65\,s &       218.03\,s &    117.71\,s \\
				2\,744 &   37\,229 & 119.86\,s &       788.58\,s & 1\,059.76\,s &   57\,308 & 204.04\,s &    1\,315.37\,s & 1\,633.51\,s \\
				4\,096 &   48\,232 & 156.47\,s &    1\,076.23\,s & 1\,432.02\,s &   77\,222 & 233.64\,s &    1\,748.85\,s & 2\,130.84\,s \\ \hline
		\end{tabular}}
		\caption{Strong scalability for the \textit{subd. problem solve time} and the \textit{coarse problem solve time} using different sparse direct linear solvers. The time includes the factorization of problem and the forward./backw. substitution during the Krylov iterations. For the whole \textit{solver time} see~\Cref{tab:mech-solver}}
		\label{tab:mech-detail}
	\end{table}
	
	
	\newpage
	\section{Conclusion} \label{sec:conclusion}
	
	In our tests, we were able to reduce the \textit{solver time} by over 80\% with the choice of the sparse solver. We should therefore take particular interest in the choice of the solver using \texttt{FROSch}. For the considered model problems, we recommend using \texttt{MUMPS} since it performed the best. We expect a good performance of \texttt{MUMPS} for other model problems as well.
	
	For the range of problems considered here, we did not face any memory issues with the direct solvers compared, and we did not specifically examine the memory usage. However, we expect differences in the range of possible subdomain and coarse problem sizes due to different memory demands of the different solvers.

	\section*{Acknowledgments}
	The authors acknowledge the DFG project 441509557 (\url{https://gepris.dfg.de/gepris/projekt/441509557}) within the Priority Program SPP2256 ``Variational Methods for Predicting Complex Phenomena in Engineering Structures and Materials'' of the Deutsche Forschungsgemeinschaft (DFG). The authors gratefully acknowledge the Gauss Centre for Supercomputing e.V. (www.gauss-centre.eu) for providing computing time on the GCS Supercomputer JUWELS~\cite{JUWELS} at the Jülich Supercomputing Centre (JSC) for the project \textit{FE2TI - High performance computational homogenization software for multi-scale problems in solid mechanics}. 
	
	\bibliographystyle{plain}
	\bibliography{all-lit.bib}
	
\end{document}